 \documentclass[final,3p,11pt,times]{elsarticle}
 
 \pdfoutput=1

\usepackage{amssymb}
\usepackage{amsmath}
\usepackage{mathtools}
\usepackage{amsthm}
\usepackage{mathrsfs}

\usepackage{array}
\usepackage{multirow}
\usepackage{pifont}
\newcommand{\cmark}{\ding{51}}%
\newcommand{\xmark}{\ding{55}}%
\usepackage{subfig}

\usepackage{float}
\usepackage{graphicx}
\usepackage{hyperref}

\numberwithin{equation}{section}
\newcommand{\mbf}[1]{\mathbf{#1}}

 \newcommand{\sgn}{\operatorname{sgn}}

\usepackage[nameinlink,noabbrev]{cleveref}
\Crefname{equation}{eq.}{eqs.} 
\Crefname{equation}{Eq.}{Eqs.}

\usepackage{color}
\PassOptionsToPackage{normalem}{ulem}
\usepackage{ulem}

\definecolor{mygreen}{rgb}{0,.6,0}

\allowdisplaybreaks

 \makeatletter
\def\ps@pprintTitle{%
  \let\@oddhead\@empty
  \let\@evenhead\@empty
  \def\@oddfoot{\reset@font\hfil\thepage\hfil}
  \let\@evenfoot\@oddfoot
}
\makeatother

\begin{document}

\begin{frontmatter}



\title{Gaussian-Like Immersed-Boundary {Kernels} with Three Continuous Derivatives and Improved Translational Invariance\tnoteref{tn1}}
\tnotetext[tn1]{This manuscript is an update to the the published article \cite{Bao2016} with a new $\mathscr{C}^3$ 5-point IB kernel. }

\author[cims]{Yuanxun Bao\corref{cor1}}
\ead{billbao@cims.nyu.edu}
\author[cims]{Alexander D.~Kaiser}
\ead{kaiser@cims.nyu.edu}
\author[cims]{Jason Kaye}
\ead{jkaye@cims.nyu.edu}
\author[cims]{Charles S.~Peskin}
\ead{peskin@cims.nyu.edu}
\address[cims]{Courant Institute of Mathematical Sciences, New York University, 251 Mercer Street, New York, NY, 10012, USA}

\cortext[cor1]{Corresponding author}

\begin{abstract}
The immersed boundary (IB) method is a general mathematical framework for studying problems involving fluid-structure interactions in which an elastic structure is immersed in a viscous incompressible fluid. In the IB formulation, the fluid described by Eulerian variables is coupled with the immersed structure described by Lagrangian variables via the use of the Dirac delta function. From a numerical standpoint, the Lagrangian force spreading and the Eulerian velocity interpolation are carried out by a regularized, compactly supported discrete delta function, which is assumed to be a tensor product of a single-variable immersed-boundary kernel. IB kernels are derived from a set of postulates designed to achieve approximate grid translational invariance, interpolation accuracy and computational efficiency. 
{In this note, we present new 5-point and 6-point immersed-boundary kernels that are $\mathscr{C}^3$ and yield a substantially improved translational invariance compared to other common IB kernels.}
\end{abstract}


\end{frontmatter}




\section{Introduction}
The immersed boundary (IB) method was proposed to study flow patterns around heart valves \cite{Peskin1977}. In the IB formulation, a viscous incompressible fluid described by Eulerian variables is assumed to occupy the entire domain, which contains an immersed structure, described by Lagrangian variables, that moves with the fluid and exerts a force on the fluid. In the spatially discretized setting, the fluid domain is represented by a uniform Eulerian grid and the immersed structure is configured as a collection of Lagrangian points or markers. The  IB kernel plays a key role in communicating between the Eulerian and Lagrangian grids by spreading applied forces to the fluid and interpolating Lagrangian marker velocity. There are three main criteria for constructing an ideal IB kernel: grid translational invariance, interpolation accuracy and computational efficiency. It is a desirable property of an IB kernel to perform force spreading and velocity interpolation that are independent of the position of Lagrangian markers relative to the Eulerian computational grid. In this case, if the IB method were applied to a translation-invariant linear system like the Stokes equations on a periodic domain,  the results would remain the same despite shifts in position of Lagrangian markers relative to the Eulerian grid \cite{Peskin2002}. There are functions that might serve as candidates for an IB kernel in terms of exact grid translational invariance; for example, the sinc function $\sin(x)/x$. However, the sinc function is not computationally efficient because its support is unbounded. In fact, as we discuss later, exact grid translational invariance is inconsistent with the assumption of compact support \cite{Peskin2002}. The process of constructing a computationally efficient IB kernel that simultaneously has good interpolation accuracy and translational invariance is non-trivial. 

In the IB method, the 3D discrete delta function is assumed to be represented by a tensor product of a single-variable kernel $\phi(r)$,
\begin{equation}
  \delta_h(\mbf{x}) = \frac{1}{h^3} \phi \left( \frac{x_1}{h} \right) \phi \left( \frac{x_2}{h} \right) \phi \left( \frac{x_3}{h} \right),
\end{equation}
where $x_1, x_2, x_3$ are the Cartesian components of $\mbf{x}$ and $h$ is the meshwidth. This representation is not essential, but it significantly simplifies the discussion, since the single-variable kernel  $\phi(r)$ is the object of interest. We first require that  $\phi(r)$ be continuous for all real $r$ and have compact support, i.e., $\phi(r) = 0$ for $|r| \geq r_s$, where $r_s$ is the radius of support. Continuity of $\phi$ is assumed in order to avoid sudden jumps in the interpolated velocity or applied force as Lagrangian markers move through the Eulerian grid. It turns out that most IB kernels are $\mathscr{C}^1$ even though the higher regularity is not explicitly assumed. The reason for that is still a mystery, but higher regularity of the IB kernel is a nice feature to have in certain applications, such as the interpolation of derivatives or the spreading of a force dipole. Compact support of $\phi$ is required for computational efficiency. 

The function $\phi(r)$ is constructed by requiring a subset of the following moment conditions:
\begin{equation*}
\displaystyle
\begin{array}{r r >{\displaystyle}l}
  \text{(i)} & \text{Zeroth moment: } & \quad \sum_{j} \phi(r-j)= 1 ,  \\ [1.5em]
  \text{(ii)} & \text{Even-odd: } & \quad \sum_{j \text{ even}} \phi(r-j) = \sum_{j \text{ odd}} \phi(r-j) = \frac{1}{2}, \\ [1.5em]
  \text{(iii)} & \text{First moment: } & \quad \sum_{j} (r-j) \, \phi(r-j) = 0, \\ [1.5em]
  \text{(iv)} & \text{Second moment: } & \quad \sum_{j} (r-j)^2 \phi(r-j) = K, \text{ for some constant } K , \\ [1.5em]
  \text{(v)} & \text{Third moment: } & \quad \sum_{j} (r-j)^3 \phi(r-j) = 0. \\ [1.5em]
\end{array}
\end{equation*}
The motivation of imposing moment conditions is well discussed in \cite{Peskin2002,Mori2012}. Briefly, the zeroth moment condition implies that the total force is the same in the Eulerian and Lagrangian grids when $\delta_h$ is used for force spreading. The even-odd condition implies (i), and was originally proposed to avoid the ``checkerboard" instability that may arise from using a collocated-grid fluid solver. Liu and Mori \cite{Mori2012} generalized this condition to the so called ``smoothing order" condition and showed that it has the effect of suppressing high-frequency errors and preventing Gibbs-type phenomena. Conservation of total torque relies on the first moment condition. Moreover, (i) and (iii) guarantee that a smooth function is interpolated with second-order accuracy when $\delta_h$ is used for interpolation. The second moment condition with $K=0$ and the third moment condition are needed to derive kernels with a higher order of interpolation accuracy. 

In addition to moment conditions, $\phi(r)$ is required to satisfy the sum-of-squares condition,
\begin{equation}
\displaystyle
 \quad \sum_{j} \left(  \phi(r-j) \right)^2 = C, \text{ for some constant } C. \label{sumofsquare}
\end{equation}
The sum-of-squares condition \Cref{sumofsquare} is a weaker version of exact grid translational invariance,
\begin{equation}
 \tilde{\phi}(r_1, r_2) = \sum_j \phi(r_1 - j) \, \phi(r_2 - j) = \Phi(r_1 - r_2), \text{ for all } r_1, r_2. \label{exact_trans_invar}
\end{equation} 
In other words, the coupling of $\phi(r)$ between any arbitrary two points $r_1, r_2$ is a function of $r_1 - r_2$ only. However, it can be shown that \Cref{exact_trans_invar} is inconsistent with the assumption of $\phi$ being compactly supported \cite{Peskin2002}. The sum-of-squares condition does give some information about the coupling of $\phi$, since it can be deduced from the Cauchy-Schwarz inequality that
\begin{equation}
  \left| \tilde{\phi}(r_1, r_2) \right| = \left| \sum_j \phi(r_1 - j) \, \phi(r_2 - j)  \right| \leq C, \text{ for all } r_1, r_2. \label{sumofsquare_ineq}
\end{equation}
\Cref{sumofsquare_ineq} guarantees that the coupling between two Lagrangian markers is strongest when the markers coincide, and furthermore \Cref{sumofsquare} implies that the self-coupling is independent of the marker position.

\begin{table}[h]
\centering
    \begin{tabular}{ | c | c | c | c | c | c | c | c | c|}
    \hline
    \multirow{2}{*}{IB Kernel}  & Support & \multirow{2}{*}{Even-Odd} & Zeroth & First  & Second  & Third & Sum of & \multirow{2}{*}{Regularity}  \\
     & $r_s$ &  & Moment & Moment & Moment & Moment & Squares & \\
    \hline
     Standard &  \multirow{2}{*}{1.5} & \multirow{2}{*}{\xmark} & \multirow{2}{*}{\cmark} & \multirow{2}{*}{ \cmark} & \multirow{2}{*}{\xmark} & \multirow{2}{*}{\xmark} & \multirow{2}{*}{$\frac{1}{2}$} & \multirow{2}{*}{{$\mathscr{C}^1$}} \\
     3-point & & & & & & & &\\
    \hline
    Smoothed & \multirow{2}{*}{2} &  \multirow{2}{*}{{\xmark}} &  \multirow{2}{*}{{\cmark}} &  \multirow{2}{*}{{\cmark}} & \multirow{2}{*}{{\xmark}} & \multirow{2}{*}{{\xmark}} & \multirow{2}{*}{{\xmark}} & \multirow{2}{*}{{$\mathscr{C}^2$}} \\
    {3-point} & & & & & & & & \\
    \hline
    Standard &  \multirow{2}{*}{2} & \multirow{2}{*}{\cmark} & \multirow{2}{*}{\cmark} & \multirow{2}{*}{ \cmark } & \multirow{2}{*}{\xmark} & \multirow{2}{*}{\xmark} & \multirow{2}{*}{$\frac{3}{8}$} & \multirow{2}{*}{{$\mathscr{C}^1$}}\\
    4-point & & & & & & & & \\
    \hline 
    {Smoothed} & \multirow{2}{*}{{2.5}} &  \multirow{2}{*}{{\xmark}} &  \multirow{2}{*}{{\cmark}} &  \multirow{2}{*}{{\cmark}} & \multirow{2}{*}{{\xmark}} & \multirow{2}{*}{{\xmark}} & \multirow{2}{*}{{\xmark}} & \multirow{2}{*}{{$\mathscr{C}^2$}} \\
    {4-point} & & & & & & & & \\
    \hline
    Standard & \multirow{2}{*}{3} & \multirow{2}{*}{\cmark} & \multirow{2}{*}{\cmark} &\multirow{2}{*}{\cmark} & \multirow{2}{*}{0} & \multirow{2}{*}{\cmark} & \multirow{2}{*}{$\frac{67}{128}$} & \multirow{2}{*}{{$\mathscr{C}^1$}} \\
    6-point & & & & & & & & \\
    \hline
     { {New }} & \multirow{2}{*}{2.5} & \multirow{2}{*}{\xmark} &  \multirow{2}{*}{\cmark} & \multirow{2}{*}{\cmark} & \multirow{2}{*}{ { $\frac{38}{60}-\frac{\sqrt{69}}{60}$ }} & \multirow{2}{*}{\cmark} & \multirow{2}{*}{ { { $\approx .393$} }} & \multirow{2}{*}{{$\mathscr{C}^3$}}  \\
     { { 5-point }}  & & & & & & & & \\ 
    \hline 
    New & \multirow{2}{*}{3} & \multirow{2}{*}{\cmark} &  \multirow{2}{*}{\cmark} & \multirow{2}{*}{\cmark} & \multirow{2}{*}{$\frac{59}{60} - \frac{\sqrt{29}}{20}$} & \multirow{2}{*}{\cmark} & \multirow{2}{*}{$\approx .326$} & \multirow{2}{*}{{$\mathscr{C}^3$}}\\
     6-point & & & & & & & & \\
    \hline 
    \end{tabular}
    \caption{  { Common immersed-boundary kernels with their properties and moment conditions they satisfy. \cmark: the kernel satisfies the moment condition; \ \xmark: the kernel does not satisfy the moment condition. In the second moment column, the value of the second moment constant $K$ is given when the second moment condition is satisfied. The regularity column shows the number of continuous derivatives each IB kernel has.}}
    \label{tab:tableIBKernels}
\end{table}

In \autoref{tab:tableIBKernels}, we list some common IB kernels and the conditions they satisfy. The most widely used IB kernel is the standard 4-point kernel \cite{Peskin2002}. The standard 3-point kernel satisfies the zeroth moment condition but not the even-odd condition. It was first introduced in an adaptive IB method using the staggered-grid discretization \cite{Roma1999}. The standard 6-point kernel (with $K=0$) satisfies all the moment conditions listed above \cite{Stockie1997}. It can be shown that the standard 6-point kernel interpolates cubic functions exactly and smooth functions with fourth-order accuracy. 
{However, as shown in \autoref{fig:test_trans_invar}, it has a larger deviation from translational invariance for a pair of markers with distance $d \approx 2.5$ even compared to the standard 4-point kernel. }
In terms of its defining postulates, our new 6-point kernel differs from the standard 6-point kernel only in the nonzero second-moment constant $K$ (the sum-of-squares constant $C$ is determined once $K$ is fixed). The new 5-point kernel assumes the same postulates as the new 6-point kernel except for  the ``even-odd'' condition. 

{The special choice of $K$ given in \Cref{specialK5pt} and \Cref{specialK} lead to new 5-point and 6-point kernels that are $\mathscr{C}^3$ and significantly improve translational invariance compared with other IB kernels.} The construction of an IB kernel with a positive and constant second moment $K$ was originally motivated by the important physical implication of the second moment in particle suspensions, namely it is associated with the quadrupole correction in the Fax\'{e}n relation for a rigid sphere in an arbitrary Stokes flow \cite{Brady1988}. The result that the new {kernels have} three continuous derivatives is unexpected, however, this makes the new kernel more generally useful.
 {By applying a smoothing technique to the standard IB kernels, Yang, {\it{et al.}} \cite{Yang2009} developed a family of $\mathscr{C}^2$ IB kernels whose first derivative satisfies up to the second moment condition for the derivative. They showed that these derivative moment conditions are intrinsically linked to the error of force spreading in the IB scheme, and IB kernels that satisfy these conditions can significantly reduce non-physical spurious oscillations of force spreading in moving-boundary problems. 
 {By differentiating the moment conditions satisfied by $\phi(r)$, we can verify that the derivative of our new $\mathscr{C}^3$ kernels satisfy up to third moment conditions as advocated by Yang {\it et al}. \cite{Yang2009}}. We will also include the smoothed 3-point and 4-point kernels in the comparison of translational invariance in \autoref{section:numerical_test}.}
 Liu and Mori developed a MATLAB routine that automatically generates all the standard IB kernels as well as many others \cite{Mori2012}. We have also made our MATLAB codes for generating the new kernels available at \url{https://github.com/stochasticHydroTools/IBMethod}.

\section{Two new kernels}

\subsection{A new 5-point kernel}
Our new 5-point kernel satisfies the sum-of-squares condition \Cref{sumofsquare} and the moment conditions (i), (iii)-(v), but it does not satisfy the even-odd condition (ii). 
The support is defined to be five grid points, i.e., $r_s = \frac{5}{2}$. 
We follow a similar derivation of the standard 4-point kernel \cite{Peskin2002}. 
By first restricting $r \in \left[-\frac{1}{2}, \frac{1}{2} \right]$, we have 5 unknowns:
\begin{align} 
\left\{ \phi(r-2) ,\,  \phi(r-1) , \,  \phi(r) , \,  \phi(r+1) , \, \phi(r+2) \right\}.
\end{align}
Note that the moment conditions (i), (iii)-(v) are four linear equations in these five unknowns, and we can express all the other four unknowns in terms of $\phi(r)$, 
\begin{align}
\phi(r-2) \ = \ & \frac{1}{12} \left(2 \phi(r)+3 K r+2 K+r^3+2 r^2-r-2\right),  \label{phi_m2} \\ 
\phi(r-1) \ = \ & \frac{1}{6} \left(-4 \phi(r)-3 K r-K-r^3-r^2+4r+4\right), \label{phi_m1} \\ 
\phi(r+1) \ = \ & \frac{1}{6} \left(-4 \phi(r)+3 K r-K+r^3-r^2-4 r+4\right), \label{phi_p1}\\ 
\phi(r+2) \ = \ & \frac{1}{12} \left(2 \phi(r)-3 K r+2 K-r^3+2r^2+r-2\right). \label{phi_p2} 
\end{align}
For the special value $r=1/2$  so that $\phi \left( \frac{5}{2} \right) = 0$, we  get
\begin{align}
\left\{ \phi \left(-\frac{3}{2}\right), \, \phi \left(-\frac{1}{2} \right), \, \phi \left(\frac{1}{2} \right), \, \phi \left(\frac{3}{2} \right) \right\} 
= \left\{\frac{1}{16} (4 K-1),\frac{1}{16} (9-4 K),\frac{1}{16} (9-4 K),\frac{1}{16} (4 K-1)\right\}. \label{phivalues}
\end{align}
Substituting these values into the sum of squares condition \Cref{sumofsquare}, we obtain an expression for $C(K)$, 
\begin{align}
C(K)  
&=\frac{1}{128} (9-4 K)^2+\frac{1}{128} (4 K-1)^2,  \label{CofK}
\end{align}
where the value of $K$ remains to be determined. 
Next we solve the quadratic equation of $\phi(r)$ from the sum-of-squares condition \Cref{sumofsquare} by using \Cref{phi_m2,phi_m1,phi_p1,phi_p2} and \Cref{CofK},
\begin{equation}
\phi(r) \ = \  \frac{1}{280}\left( -40 K-40 r^2+136  + \sqrt{ 2\beta(r) + 2\gamma(r)} \, \right), \label{phideffivepoint} 
\end{equation}
where 
\begin{align}
\beta(r) & \ = \   -12600 K^2 r^2+3600 K^2-8400 K r^4+25680 K r^2-6840 K+3123, \\ 
\gamma(r) & \ = \   -40 r^2 \left(35 r^4-202 r^2+311\right).
\end{align}
We note that the positive square root is chosen in \Cref{phideffivepoint} because of the continuity assumption of $\phi$, i.e., by setting $r= \frac{1}{2}$, we select the branch that gives $\phi \left(\frac{1}{2} \right) = \frac{1}{16}(9-4K)$ as in \Cref{phivalues}.  

We have the freedom to choose $K$ to construct $\phi$ with higher regularity. 
A symbolic calculation in \emph{Mathematica} matching derivatives of $\phi(r)$ at $r = \frac{1}{2}$ reveals that for any $K \in \left[ 0, \frac{21}{20} \right)$, $\phi \in \mathscr{C}^1$. 
For the second derivative to be continuous at $r = \frac{1}{2}$, $K$ must satisfy 
\begin{align} 
720 K^2-912 K+275 = 0. \label{secondderivative_condition}
\end{align} 
We can verify (by plotting) that exactly one of the two roots of \Cref{secondderivative_condition}
guarantees that $\beta(r) + \gamma(r) \geq 0$,  so that $\phi(r)$ is  real-valued for  $r \in \left[ \frac{1}{2}, \frac{1}{2}\right]$, and this special value of $K$ is 
\begin{align}
K =  \frac{1}{60} \left(38-\sqrt{69} \right).  \label{specialK5pt}
\end{align}
If we  proceed further with matching the third derivative to be continuous at $r = \frac{1}{2}$, then $K$ must satisfy 
\begin{align} 
60480 K^3-116208 K^2+73260 K-15125 = 0, \label{roots_secondderivative}
\end{align} 
whose roots are
\begin{align} 
K = \left \{ \frac{55}{84}, \ \frac{1}{60} \left(38 + \sqrt{69} \right), \ \frac{1}{60} \left(38 - \sqrt{69}  \right) \right\}. \label{roots_thirdderivative} 
\end{align} 
We observe that our choice of $K$ in \Cref{specialK5pt} is among the roots of the third derivative matching condition, and therefore, it also makes the third derivative of $\phi$ continuous at $r = \frac{1}{2}$.

Remarkably, the derivative matching conditions at $r = \frac{1}{2}$ are sufficient to ensure that $\phi \in \mathscr{C}^{3}$ everywhere. 
Existence and continuity of derivatives was never assumed a priori, but was only a consequence of an appropriate choice of $K$. 
Moreover, note that for $K \in \left[ 0, \frac{1}{60} \left(38-\sqrt{69} \right) \right)$, we have  $\phi''\left(\frac{5}{2}\right) < 0$. 
Since $\phi \left(\frac{5}{2}\right) = \phi' \left(\frac{5}{2}\right) = 0$, this implies that $\phi(r)$ takes negative values in a neighborhood of $r = \frac{5}{2}$.  We emphasize that the special choice of $K$ given by \Cref{specialK5pt} is the smallest positive $K$ for which $\phi$ is non-negative, and it is also the only value of $K$ that gives three continuous derivatives of $\phi$. 

\subsection
{A new 6-point kernel}

\begin{figure}
     \centering
     \subfloat[][]{\includegraphics[width= .45\linewidth]{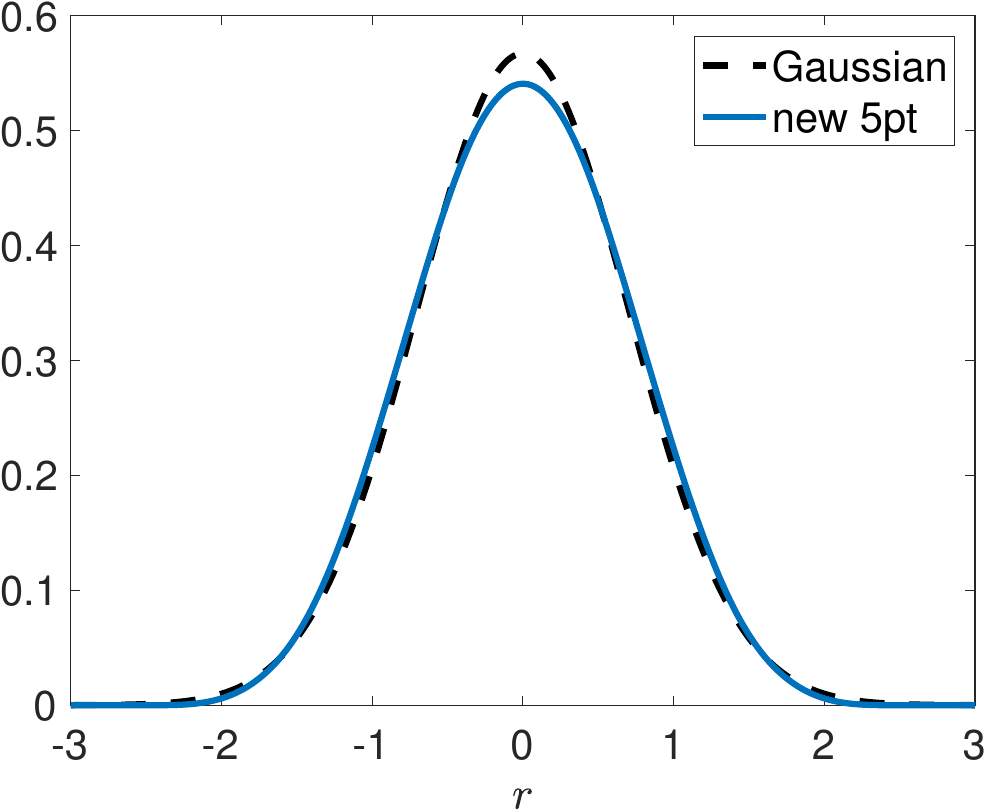}}  \hspace{1em}
     \subfloat[][\label{fig:new5pt_derivatives}]{\includegraphics[width= .44\linewidth]{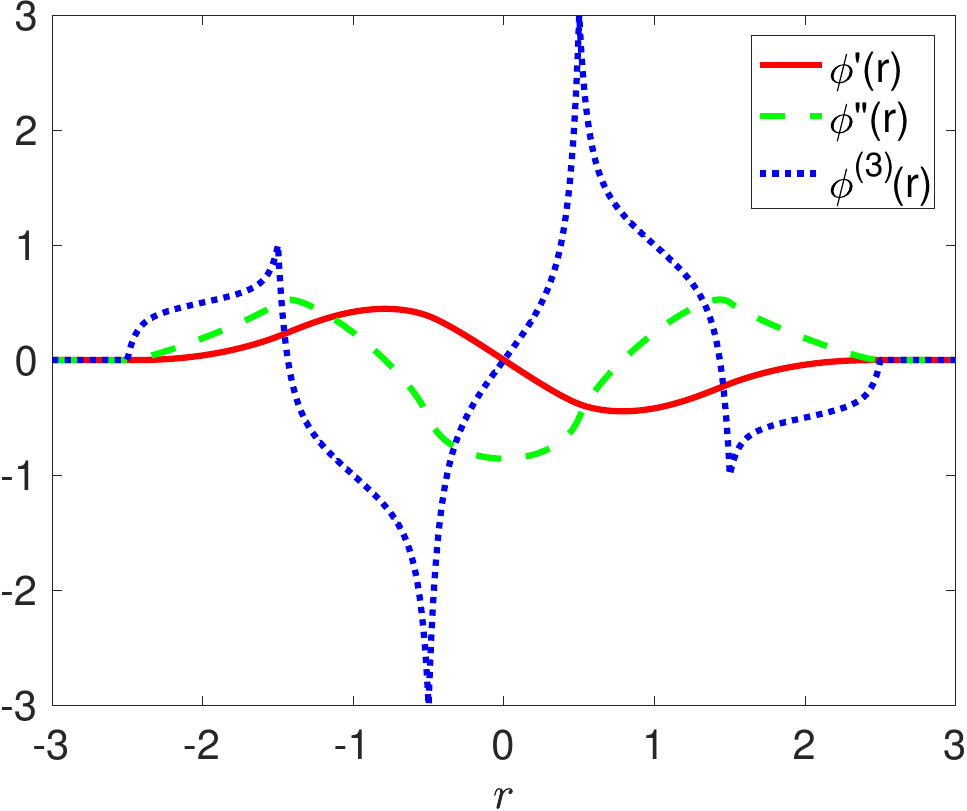}}
     \caption{(a) The new 5-point kernel compared with the Gaussian with the second moment given by \Cref{specialK5pt}. (b) The first three derivatives of the new 5-point kernel.}
     \label{fig:new5pt_plot} 
\end{figure}

\begin{figure}
     \centering
     \subfloat[][]{\includegraphics[width= .45\linewidth]{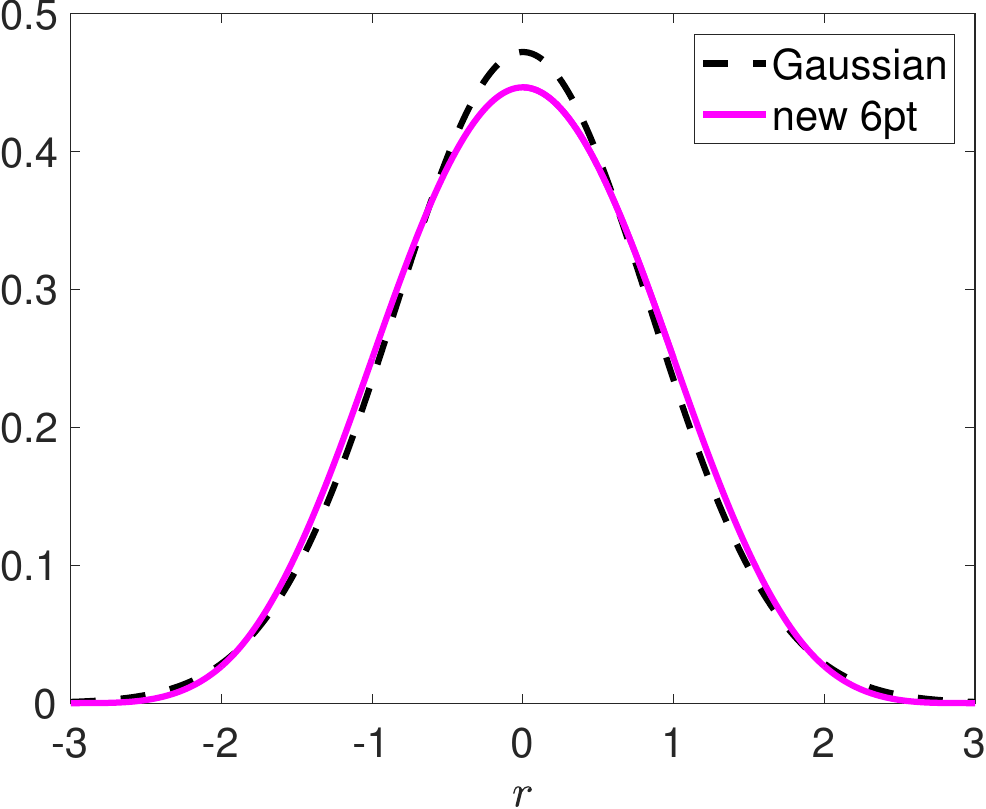}}  \hspace{1em}
     \subfloat[][\label{fig:new6pt_derivatives}]{\includegraphics[width= .46\linewidth]{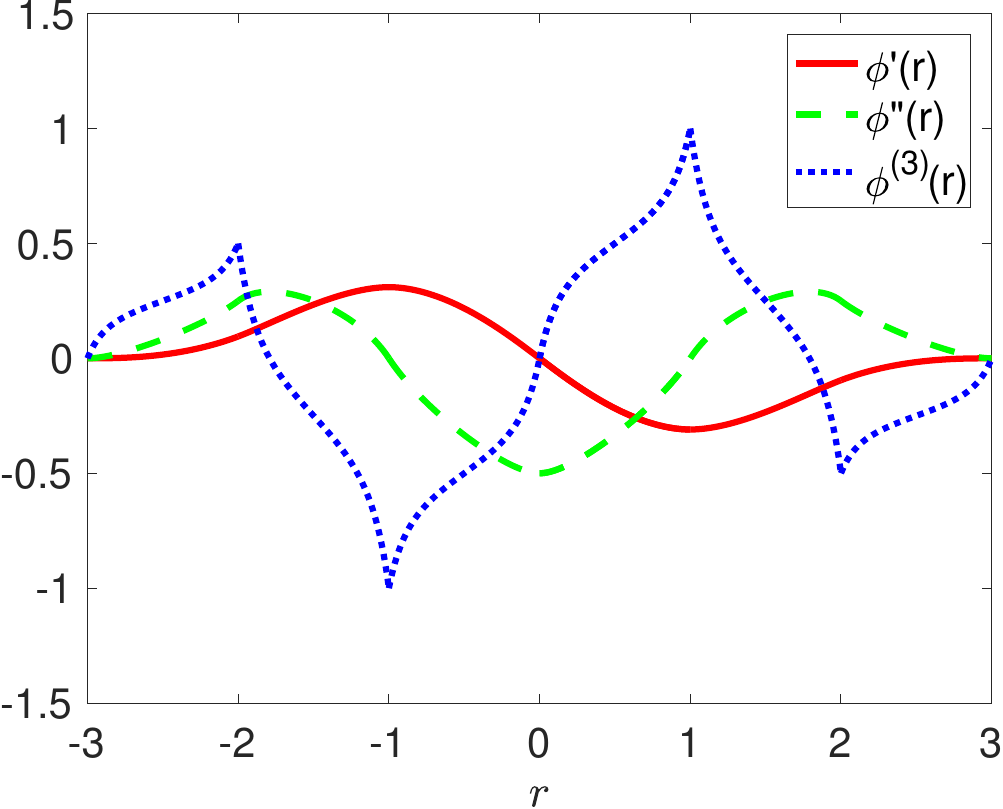}}
     \caption{(a) The new 6-point kernel compared with the Gaussian with the second moment given by \Cref{specialK}. (b) The first three derivatives of the new 6-point kernel.}
     \label{fig:new6pt_plot} 
\end{figure}

\begin{figure}
     \centering
     \subfloat[][]{\includegraphics[width= .45\linewidth]{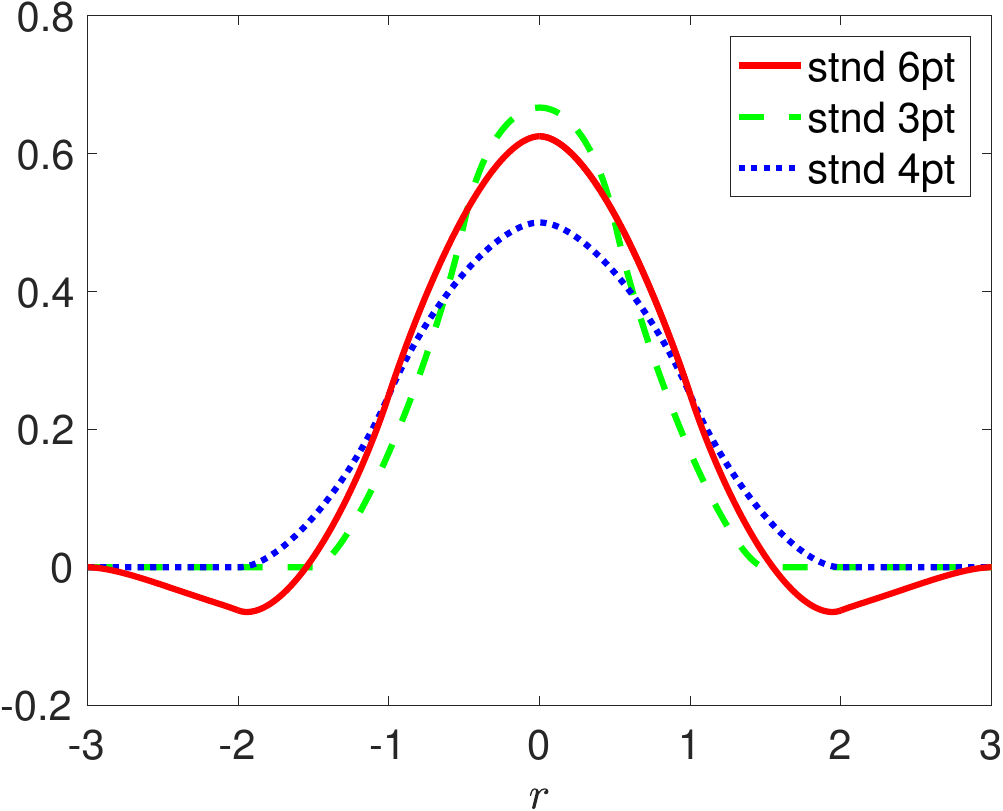}\label{fig:stdkernels}} 
     \subfloat[][]{\includegraphics[width= .45\linewidth]{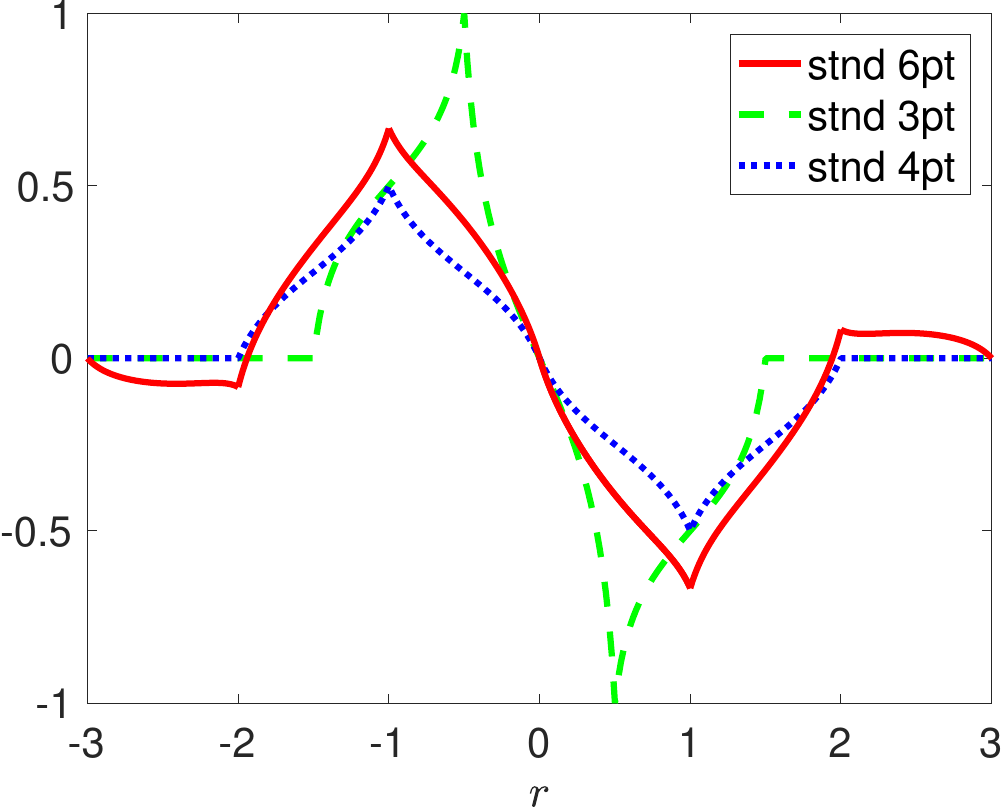}\label{fig:stdkernels_d}}
     \caption{(a) The standard 3-point, 4-point, and 6-point kernels. (b) The first derivatives of the standard 3-point, 4-point, and 6-point kernels.}
     \label{fig:stdkernels_fig}
\end{figure}

Our new 6-point kernel satisfies the sum-of-squares condition \Cref{sumofsquare} and the moment conditions (ii)-(v) (and therefore (i)) with the second-moment constant 
\begin{equation}
 K = \frac{59}{60} - \frac{\sqrt{29}}{20}. \label{specialK}
\end{equation}
The derivation of this kernel follows a nearly identical procedure to that of the new 5-point kernel.
First, the sum-of-squares constant $C$ can be expressed in terms of $K$ by considering the special case $r = 0$. Next, by restricting $r$ to the interval $[0,1]$, we have 6 unknowns: $\phi(r-3), \phi(r-2)$, $\phi(r-1), \phi(r), \phi(r+1), \phi(r+2)$ and 6 equations (the even-odd condition accounts for two equations). By expressing all the other unknowns in terms of $\phi(r-3)$ using (ii)-(v), we can solve for $\phi(r-3)$ from the quadratic equation determined by \Cref{sumofsquare}. The continuity assumption of $\phi$ is now used to select the appropriate root to piece together a continuous function, i.e., by setting $r=0$, we select the root that gives $\phi(-3) = 0$. As mentioned earlier, $\phi$ being $\mathscr{C}^1$ follows implicitly from our defining postulates, i.e., $\phi'(-3) = 0$.  We have the freedom to choose $K$ so that $\phi''(-3) = 0$, which uniquely determines the special value  \Cref{specialK}. 
The formula for the new 6-point kernel is given by
\begin{align}
 \beta(r)   \ = \ & \frac{9}{4} - \frac{3}{2}(K+r^2) + ( \frac{22}{3}-7K )r - \frac{7}{3}r^3 \, , \\ 
 \gamma(r)  \ = \ & -\frac{11}{32}r^2 + \frac{3}{32}(2K+r^2)r^2 + \nonumber \\ 
                & \frac{1}{72}\left((3K-1)r+r^3\right)^2  + \frac{1}{18} \left((4-3K)r - r^3\right)^2 \, ,
\end{align}
\begin{align}
\phi(r-3)  \ = \  &  \frac{-\beta(r) + \sgn\left(\frac{3}{2}-K\right) \sqrt{\beta^2(r) - {112} \gamma(r)} }{56} \, ,\\
 \phi(r-2)  \ = \ & -3\phi(r-3)  -  \frac{1}{16}  +  \frac{K+r^2}{8} + \frac{(3K-1)r}{12}  +  \frac{r^3}{12} \, ,  \\ 
 \phi(r-1)  \ = \ & 2 \phi(r-3)  +  \frac{1}{4}   +  \frac{(4-3K)r}{6} - \frac{r^3}{6} \, , \\ 
 \phi(r)    \ = \ & 2 \phi(r-3)  +  \frac{5}{8}   -  \frac{K+r^2}{4} \, , \\ 
 \phi(r+1)  \ = \ & -3\phi(r-3)  +  \frac{1}{4}   -  \frac{(4-3K)r}{6} + \frac{r^3}{6} \, , \\ 
 \phi(r+2)  \ = \ & \phi(r-3)    -  \frac{1}{16}  +  \frac{K+r^2}{8} - \frac{(3K-1)r}{12}  - \frac{r^3}{12} \, .
\end{align}
Note that, in the formula presented above, $r \in [0,1]$.
The new 5-point and 6-point kernels are Gaussian-like function, as shown in  \autoref{fig:new5pt_plot} and \autoref{fig:new6pt_plot}, and they both have three continuous derivatives. As a comparison, the standard 3-point, 4-point, 6-point kernels and their continuous first derivative are plotted in \autoref{fig:stdkernels_fig}. 
We notice that the new 5-point and 6-point kernels are non-negative for all $r$, whereas the standard 6-point kernel has negative tails.

%
%
\begin{figure}
     \centering
     \subfloat[][\label{fig:transinvar1}]{\includegraphics[width= .5\linewidth]{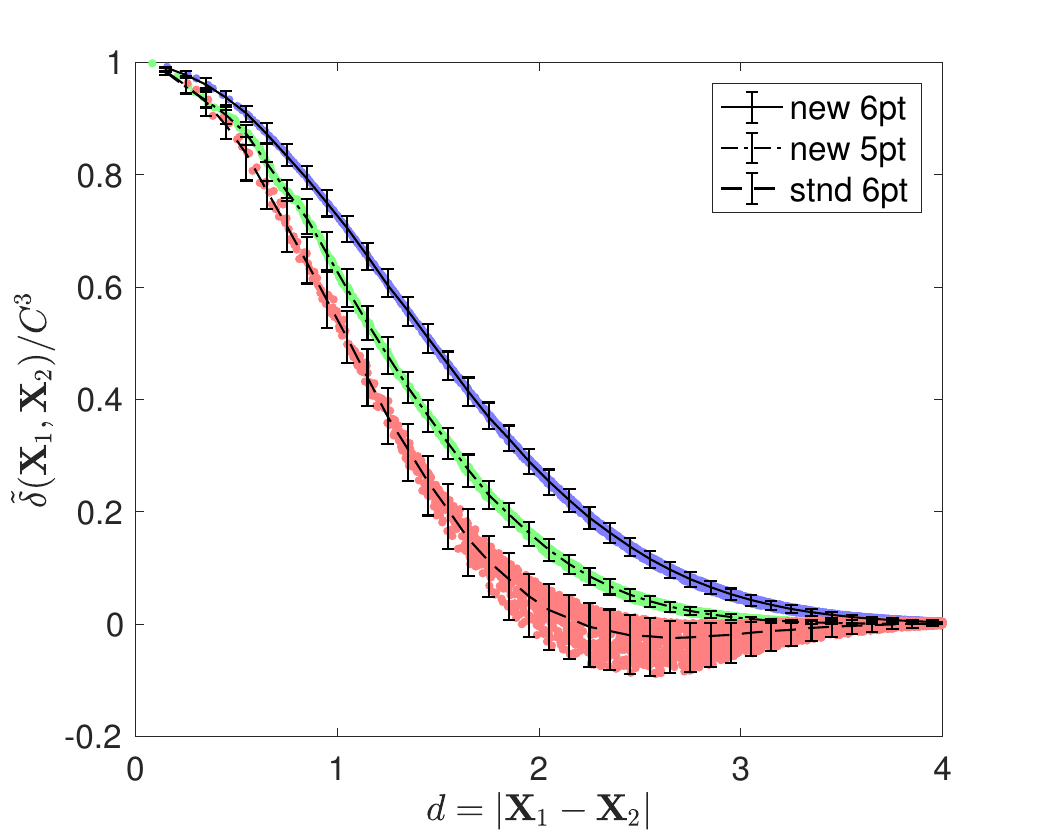}}
     \subfloat[][\label{fig:transinvar2}]{\includegraphics[width= .5\linewidth]{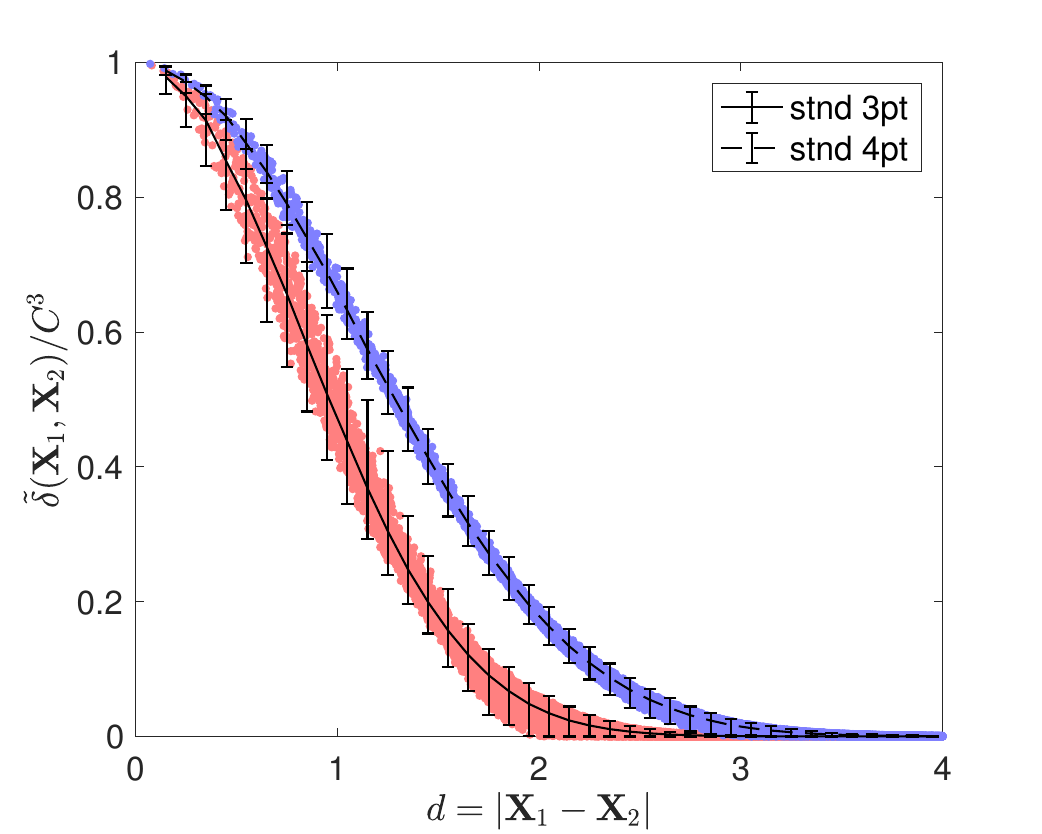}}\\
     \subfloat[][\label{fig:transinvar3}]{\includegraphics[width= .5\linewidth]{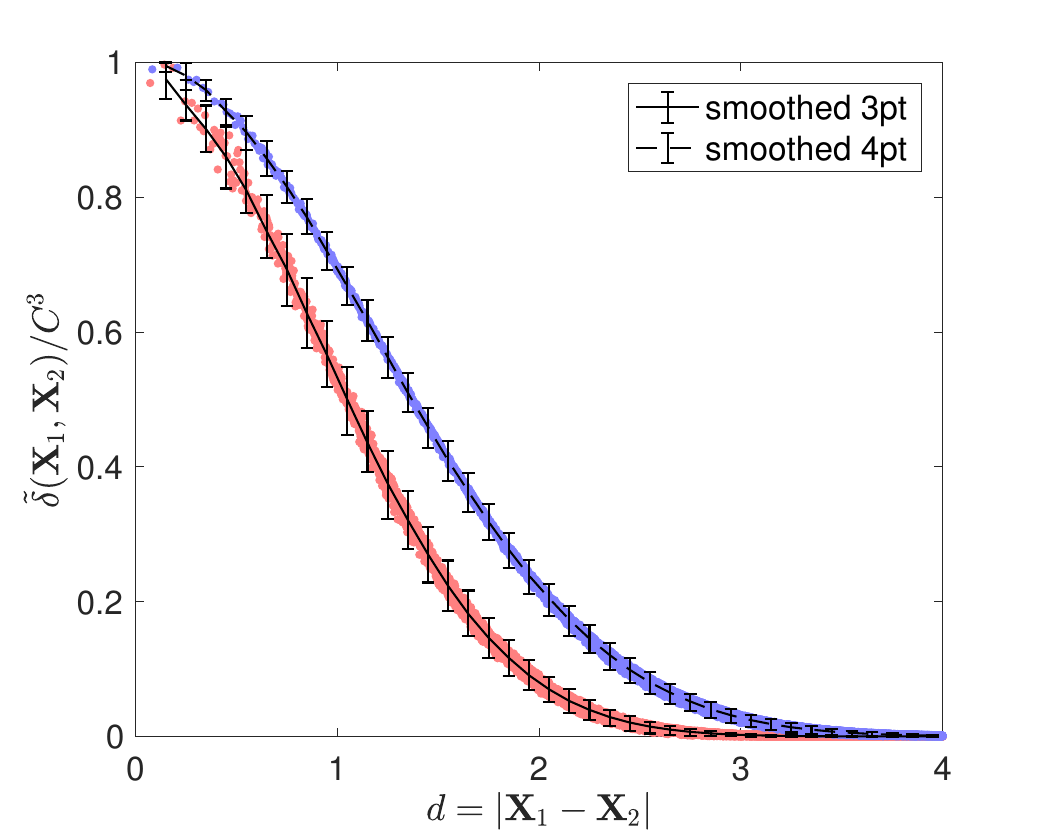}}
     \caption{Normalized $\tilde{\delta}(\mbf{X}_1, \mbf{X}_2) $ is plotted versus $d=|\mathbf{X}_1 - \mathbf{X}_2|$ for $10^5$ pairs of randomly selected Lagrangian markers. The data are binned according to $d=|\mbf{X}_1 - \mbf{X}_2|$ and error-bars showing the maximum, mean and minimum of each bin are overlaid with the data. The deviation in the data gives a quantitative measure of translational invariance of an IB kernel. The standard deviation in the data for the new 6-point kernel (blue) is an order of magnitude smaller than for the standard IB kernels. (a) The standard 6-point kernel vs. the new 5-point kernel vs. the new 6-point kernel. (b) The standard 3-point kernel vs. the standard 4-point kernel. {(c) The smoothed 3-point kernel vs. the smoothed 4-point kernel. }}
     \label{fig:test_trans_invar}
\end{figure}


\section{Numerical Tests} \label{section:numerical_test}
In this section, we demonstrate a significant improvement in the translational invariance\footnote{{The test that we use actually checks for rotational invariance at the same time, since it involves the Euclidean distance between a pair of markers, and not merely the vector from one marker to the other.} } of our new 6-point kernel. We randomly select $10^5$ pairs of Lagrangian markers $\mbf{X}_1, \mbf{X}_2$ in a periodic box $[0,32]^3$ with meshwidth $h = 1$ and compute the 3D generalization of \Cref{exact_trans_invar}, 
\begin{equation}
	\tilde{\delta}(\mbf{X}_1, \mbf{X}_2) = \sum_{\mbf{x} \in g_h } \delta_h (\mbf{x} - \mbf{X}_1) \,\delta_h (\mbf{x} - \mbf{X}_2), \label{deltaX1X2}
\end{equation}
where $\mbf{x}$ denotes a grid point on the Eulerian grid $g_h$. In \autoref{fig:test_trans_invar}, we plot $\tilde{\delta}(\mbf{X}_1, \mbf{X}_2)$ normalized by the constant $C^3$ from \Cref{sumofsquare}, versus the distance $d=\left| \mbf{X}_1 - \mbf{X}_2 \right|$. The data are binned according to $d=\left| \mbf{X}_1 - \mbf{X}_2 \right|$, and error-bars showing the maximum, mean and minimum of each bin are overlaid with the data. If an IB kernel were exactly translation-invariant,  the plot of $\tilde{\delta}(\mbf{X}_1, \mbf{X}_2)$ would be a curve. The spreading pattern in the data around this curve clearly indicates that none of the IB kernels compared here are exactly translation-invariant, but gives a qualitative picture of how close to translational invariance each kernel is. The data of the new IB kernels and {the smoothed 4-point kernel} almost form a curve, while the data of the other kernels have larger deviations from the mean. Moreover, the deviation in the data of the new IB kernels is uniform for all distances, but the standard 6-point kernel has a much larger deviation for $d \approx 2.5$. For a more quantitative comparison, we summarize the maximum standard deviation of all bins for each IB kernel in \autoref{tab:maxstd}. The maximum standard deviation of the new IB kernels is an order of magnitude smaller than that of the {standard IB} kernels, {and is about half of the deviation  of the smoothed 4-point kernel.} 

{In terms of computational cost, we summarize the computation time of $\tilde{\delta}(\mbf{X}_1, \mbf{X}_2)$ for $10^5$ pairs of Lagrangian markers using the kernels we have compared. The timings are based on simulations performed on a desktop with Intel Core i7-4770 CPU 3.40GHz under the MATLAB environment. The main cost of using an IB kernel in spreading/interpolation depends on its support size. In \autoref{tab:maxstd}, the new 5-point kernel is about two times more expensive, and the new 6-point kernel is about three-to-four times more expensive than a  4-point kernel in our comparison, because the new 5-point and 6-point kernels communicate with 125 and 216 nearby grid points respectively in spreading/interpolation in 3D, while a 4-point kernel only communicates with 64 nearby grid points. The smoothed 3-point and 4-point kernels are more expensive than their standard counterparts in that they have wider supports as shown in \autoref{tab:tableIBKernels}. In all respects, the new IB kernels achieve significant improvement in grid translational invariance with a modest increase in computational cost.}

\begin{table}[h]
\centering
\begin{tabular}{c c c c c c c c}
\hline\hline
\multirow{2}{*}{} & Standard & {Smoothed} & Standard & {Smoothed} & Standard & New & New \\ 
& 3-point & {3-point} & 4-point & {4-point} &  6-point & 5-point &  6-point \\ 
\hline
maximum & \multirow{2}{*}{0.0428} &  {\multirow{2}{*}{0.0212}} &  \multirow{2}{*}{0.0168} &  {\multirow{2}{*}{0.0083}} & \multirow{2}{*}{0.0296} & \multirow{2}{*}{0.0051} & \multirow{2}{*}{0.0042}  \\
std. dev. & & & & & &   \\
\hline
{computation} & {\multirow{2}{*}{6.34s}} & {\multirow{2}{*}{9.52s}} & {\multirow{2}{*}{9.01s}} & {\multirow{2}{*}{12.06s}} & {\multirow{2}{*}{31.19s}} & {\multirow{2}{*}{17.54s}} & {\multirow{2}{*}{30.86s}} \\
{time} & & & & & \\
\hline\hline
\end{tabular}
\caption{Maximum standard deviation of $\tilde{\delta}(\mbf{X}_1, \mbf{X}_2)$ over all bins for various IB kernels, and the computation time for computing $\tilde{\delta}(\mbf{X}_1, \mbf{X}_2)$ for $10^5$ markers. }
\label{tab:maxstd}
\end{table}


\section{Conclusion}
In this note, we have presented  new immersed-boundary kernels that are used for force spreading and velocity interpolation in the immersed boundary method. The new IB kernels are distinguished from other existing IB kernels by their nonzero second-moment constant $K$. A special choice of $K$  leads to a 5-point  or  6-point IB kernel that is $\mathscr{C}^3$ and features substantially improved translational invariance compared with the existing standard IB kernels.  {Recently, we have successfully applied the new IB kernels to a new IB method with an exactly divergence-free interpolated velocity field \cite{DivFreeIB}, in which derivatives of the discrete delta function are involved, and have achieved $10^3$ to $10^5$ times improvements in volume conservation of the IB method.} We believe that, {in many other applications in which derivatives of the IB kernel are needed}, the improved grid invariance and regularity of the new IB kernels will be worth its extra computational cost.

\section{Acknowlegements}
We thank Aleksandar Donev for many enlightening discussions on this work, and in particular for his suggestion to use the nonzero second moment condition as a postulate of the new IB kernels. Y. Bao was supported in part by the Air Force Office of Scientific Research under grant number FA9550-12-1-0356, as well as the U.S. 
Department of Energy Office of Science, Office of Advanced Scientific 
Computing Research, Applied Mathematics program under Award Number 
DE-SC0008271. 
A. Kaiser was supported by the National Science Foundation Graduate Research Fellowship Program under grant DGE-1342536. 
J. Kaye was supported in part by the National Science 
Foundation under grants NSF DMS-1115341 and DMS-1016554.





\section*{References}
\bibliographystyle{elsarticle-num}
\bibliography{IBKernelNote}

\begin{thebibliography}{1}
\expandafter\ifx\csname url\endcsname\relax
  \def\url#1{\texttt{#1}}\fi
\expandafter\ifx\csname urlprefix\endcsname\relax\def\urlprefix{URL }\fi
\expandafter\ifx\csname href\endcsname\relax
  \def\href#1#2{#2} \def\path#1{#1}\fi

\bibitem{Bao2016}
Y.~Bao, J.~Kaye, C.~S. Peskin, A gaussian-like immersed-boundary kernel with
  three continuous derivatives and improved translational invariance, Journal
  of Computational Physics 316 (2016) 139 -- 144.

\bibitem{Peskin1977}
C.~S. Peskin, Numerical analysis of blood flow in the heart, J. Computational
  Phys. 25~(3) (1977) 220--252.

\bibitem{Peskin2002}
C.~S. Peskin, The immersed boundary method, Acta Numerica 11 (2002) 479--517.

\bibitem{Mori2012}
Y.~Liu, Y.~Mori, Properties of discrete delta functions and local convergence
  of the immersed boundary method, SIAM J. Numer. Anal. 50~(6) (2012)
  2986--3015.

\bibitem{Roma1999}
A.~M. Roma, C.~S. Peskin, M.~J. Berger, An adaptive version of the immersed
  boundary method, Journal of Computational Physics 153~(2) (1999) 509 -- 534.

\bibitem{Stockie1997}
J.~M. Stockie, Analysis and computation of immersed boundaries, with
  application to pulp fibres, ProQuest LLC, Ann Arbor, MI, 1997, thesis
  (Ph.D.)--The University of British Columbia (Canada).

\bibitem{Brady1988}
J.~F. Brady, R.~J. Phillips, J.~C. Lester, G.~Bossis, Dynamic simulation of
  hydrodynamically interacting suspensions, Journal of Fluid Mechanics 195
  (1988) 257--280.

\bibitem{Yang2009}
X.~Yang, X.~Zhang, Z.~Li, G.-W. He, A smoothing technique for discrete delta
  functions with application to immersed boundary method in moving boundary
  simulations, Journal of Computational Physics 228~(20) (2009) 7821 -- 7836.

\bibitem{DivFreeIB}
Y.~X. Bao, C.~S. Peskin, B.~Griffith, D.~McQueen, A.~Donev, {An Immersed
  Boundary Method with Divergence-Free Velocity Interpolation}, submitted to J.
  Comp. Phys., ArXiv:1701.07169 (2017).

\end{thebibliography}

\end{document}